\documentclass [12pt]{article}
\usepackage{amsmath,amssymb,amscd}
\usepackage{hyperref}
\usepackage{graphicx}

\begin{document}
\def\l{\lambda}
\def\m{\mu}
\def\a{\alpha}
\def\b{\beta}
\def\g{\gamma}
\def\G{\Gamma}
\def\d{\delta}
\def\e{\epsilon}
\def\o{\omega}
\def\O{\Omega}
\def\v{\varphi}
\def\t{\theta}
\def\r{\rho}
\def\bs{$\blacksquare$}
\def\bp{\begin{proposition}}
\def\ep{\end{proposition}}
\def\bt{\begin{theo}}
\def\et{\end{theo}}
\def\be{\begin{equation}}
\def\ee{\end{equation}}
\def\bl{\begin{lemma}}
\def\el{\end{lemma}}
\def\bc{\begin{corollary}}
\def\ec{\end{corollary}}
\def\pr{\noindent{\bf Proof: }}
\def\note{\noindent{\bf Note. }}
\def\bd{\begin{definition}}
\def\ed{\end{definition}}
\def\C{{\mathbb C}}
\def\P{{\mathbb P}}
\def\Z{{\mathbb Z}}
\def\d{{\rm d}}
\def\deg{{\rm deg\,}}
\def\deg{{\rm deg\,}}
\def\arg{{\rm arg\,}}
\def\min{{\rm min\,}}
\def\max{{\rm max\,}}
 \def\d{{\rm d}}
\def\deg{{\rm deg\,}}
\def\deg{{\rm deg\,}}
\def\arg{{\rm arg\,}}
\def\min{{\rm min\,}}
\def\max{{\rm max\,}}

\newtheorem{theo}{Theorem}[section]
\newtheorem{lemma}{Lemma}[section]
\newtheorem{definition}{Definition}[section]
\newtheorem{corollary}{Corollary}[section]
\newtheorem{proposition}{Proposition}[section]

\newcommand{\epfv}{\hspace{1em}$\Box$\vspace{1em}}

\begin{titlepage}
\begin{center}

{\LARGE{\bf {Reconstruction of Planar Domains

\bigskip

from Partial Integral Measurements}}}

\vskip 6mm

D. Batenkov$^1$, V. Golubyatnikov $^{2,3}$, Y. Yomdin $^1$

\medskip

 1. Weizmann Institute of Science, Rehovot, Israel.

 2. Sobolev Institute of Mathematics, Novosibirsk, Russia.

 3. Novosibirsk State University, Russia.

\end{center}

\vspace{2 mm}
\begin{center}

{ \bf Abstract}
\end{center}
{\small{We consider the problem of reconstruction of planar domains
 from their moments. Specifically, we consider domains with boundary
 which can be represented by
 a union of a finite number of pieces whose graphs are solutions of a linear
 differential equation with polynomial coefficients. This includes domains with
 piecewise-algebraic and, in particular, piecewise-polynomial boundaries.
 Our approach is based on one-dimensional reconstruction method of \cite{bat}
 and a kind of ``separation of variables" which reduces the planar problem to
 two one-dimensional problems, one of them parametric. Several explicit
 examples of reconstruction are given.

Another main topic of the paper concerns ``invisible sets" for various types of
incomplete moment measurements. We suggest a certain point of
view which stresses remarkable similarity between several apparently unrelated
problems. In particular, we discuss zero quadrature domains (invisible for
harmonic polynomials), invisibility for powers of a given polynomial, and
invisibility for complex moments (Wermer's theorem and further developments).
The common property we would like to stress is a ``rigidity" and
symmetry of the invisible objects.

\vspace{3 mm}
\begin{center}
------------------------------------------------
\vspace{2 mm}
\end{center}
This research was supported by the ISF, Grants No. 639/09,
and by the Minerva Foundation.}}

\end{titlepage}
\newpage


\section{Introduction}
\setcounter{equation}{0}

In this paper we continue our study of nonlinear problems of reconstruction
of multidimensional objects from incomplete collection of integral measurements.
The paper has two parts, closely related but different in their goals. In the
first part we present a method of reconstruction of planar domains of a certain
special class from finite collections of moments. In the second part we discuss
the structure of sets and functions ``invisible" for a certain collection of
moment measurements.

\medskip

 In more detail, the object we would like to reconstruct is a 2-dimensional finite
 domain $ G \subset  {\mathbb R}^2$ which we assume to belong to a certain
 finite-dimen\-sio\-nal family $G_\lambda$ specified by a finite number of discrete
 and continuous parameters $\lambda$. Specifically, we shall assume that the
 boundary of $G$ is a union of a finite number of pieces whose graphs are
 solutions of a linear differential equation with polynomial coefficients.

The measurements are represented by finite collections of the moments
$m_{\alpha, \beta} $ of the characteristic function $ \chi_G (x,y)$ of the
 domain $G$:
 \be
 m_{\alpha, \beta} = \iint_{R^2} \, \chi_G (x,y)\cdot
 x^{\alpha} y^{\beta} dx dy.
 \ee
 Our main problem is to provide an explicit (and potentially efficient)
 reconstruction method, and, in particular, to estimate
 a minimal possible set
 of these moments sufficient for unique reconstruction of the domain $G$.

 Similar inverse problems have been intensively studied, including
 reconstruction from their moments of polygons, of quadrature domains, of
 certain ``dynamic" semi-algebraic sets,
  see \cite{bgy1,bgy2,ggmpv,ghmp,pp}
 and references therein. In a more general context the problem of domain
 reconstruction from its moments appears as a part of broad field of inverse
 problems in Potential Theory (see, for example, \cite{Var.Eti}). Rather
 similar questions arise in reconstruction from tomography measurements
 (\cite{GA,G,PGK}).

\medskip

 Our approach is based on one-dimensional reconstruction method of
 \cite{bat}
 applicable to piecewise continuous functions satisfying on each
 continuity interval
 a linear differential equation with polynomial coefficients. Then we use
 a kind of ``separation of variables" which reduces the planar problem to
 two one-dimensional problems, one of them parametric.

We expect that a reconstruction method for piecewise-smooth functions
given in \cite{by} can be extended in a similar way also to planar and higher
dimensional piecewise-smooth functions.

\medskip

 The second part of the present paper is devoted
 to ``invisible sets" for
 various types of incomplete measurements.
 Here we do not provide new results (besides a couple of examples),
 but rather suggest a certain point of view which stresses
 remarkable similarity between several apparently unrelated
 ``moment vanishing" problems. In particular, we discuss zero quadrature domains
 (invisible for harmonic polynomials), invisibility for polynomials annihilating
 other partial differential operators, invisibility for powers of a given
 polynomial, and invisibility for complex moments (Wermer's theorem and further
 developments). In all these cases we stress a common property of ``rigidity" and
 symmetry of the invisible objects.

\section{One-dimensional case}
\setcounter{equation}{0}

Reconstruction problem in dimension one has been settled in a pretty
satisfactory way for many important finite-dimensional families of functions.
This includes linear combinations of shifts of known functions, signals
with ``finite rate of innovation", piecewise D-finite functions which we use
below, piecewise-smooth functions, and many other cases
(see \cite{bat,bsy,by,vet,Sar.Yom} and references therein).

\subsection{Piecewise D-finite reconstruction}\label{sec:picewise-dfin-onedim}

Let $g(x)$ be a function with a support $ [a,b] \subset {\mathbb R}^1 $,
satisfying the following condition: there exists a finite set of ${\cal K}+2$
points $a = \xi_0 < \xi_1 < \ldots < \xi_{{\cal K}+1} = b$, such that on
each segment $[ \xi_n , \xi_{n+1}], \ n = 0,1, \ldots {\cal K}$,
the function $ g(x)$ is continuous and satisfies there a linear differential
equation
 \be
 D_n g(x) = \sum\limits_{j=0}^{N} p_{n,j}(x)\left( \frac{d^j g}{dx^j}\right) = 0
 \ee
 with polynomial coefficients
 $p_{n,j}(x)= \sum\limits_{i=0}^{k_{n,j}} a_{n,i,j} x^i $,
 $\, p_{N,j}\ne 0$ on $ [ a,b ]$.
 At the points $ \xi_n $ the function $ \, g(x) \, $ may have jumps. Such
 functions are described by a finite collection of discrete and continuous
 parameters, and they are called piecewise $D$-finite. Without loss of generality
(at least theoretically) we can assume that all the operators $D_n=D$ are the same.
 In particular, piecewise-algebraic, and,
 specifically, piecewise-polynomial
 functions belong to this class.
 In the last case the differential operator is
 $D= \frac{d^{N+1}}{dx^{N+1}} $,
 where $N$
 is the maximal degree
 of the polynomial pieces of $g$.

It was shown in \cite{bat} that the collection of ``discrete" parameters
$ {\cal K} $, $ N $, $\{k_{n,j} \} $, together with a sufficiently large
collection of the moments
 $$
 m_{\alpha} = \int\nolimits_{a}^{b}
 g(x)\cdot x^{\alpha} dx, \quad \alpha = 0, 1, 2, \ldots, \mu
 $$
determine uniquely any $ D $-finite function $ g(x) $ with all
the points $ \xi_0, \ldots \xi_{{\cal K}+1} $ of its possible discontinuity,
as well as the coefficients of the differential operator $ D$.

For piecewise-polynomial functions the number $\mu$ of the moments required
for reconstruction, depends only on the discrete data: the number of jumps
 ${\cal K}$ and on the maximal degree
 $N$
 of the pieces. It is shown in \cite{bat} that in piecewise-polynomial case
 \[
 \mu=\mu({\cal K},N)=
 \max\biggl\{ 2(N + 1)\mathcal{K}-2, (\mathcal{K} + 1)(N+1) \biggr\}.
 \]

 A similar, but
 more complicated expression for $\mu$ can be written in piecewise-algebraic case.
 However, for general $D$-finite functions, with respect to a general second
 (and higher) order differential operators $D$, the number $\mu$ may depend
 also on specific coefficients of $D$.

Let us give a very simple example of this latter phenomenon. Let $L_n(x)$
be the $n$-th Legendre polynomial, defined as
$L_n(x)={1\over {2^n n!}}{{d^n}\over {dx^n}}[(x^2-1)^n].$ Legendre polynomials
are pairwise orthogonal on $[-1,1]$ and they satisfy the second order Legendre
differential equation ${d\over {dx}}[(1-x^2){d\over {dx}}L_n(x)]+n(n+1)L_n(x)=0.$
Since $L_j(x), \ j\leq n,$ form a basis of the space of all polynomials of degree
$n$ we conclude that $L_n(x)$ is orthogonal to $1,x,x^2,\dots,x^{n-1}$. Hence the
 moments $m_j(L_n)=\int_{-1}^1 x^j L_n(x)dx$ vanish for $j=0,1,\dots,n-1$. We
 conclude that a $D$-finite function $L_n(x)$ on $[-1,1]$ cannot be reconstructed
 from less than $n+1$ its moments. So $\mu$ above depends not only on the order
 and degree of the Legendre operator, but also on a specific value of the
 parameter $n$ in it.

\smallskip

 Notice that the leading coefficient of the Legendre equation vanishes at both
 the endpoints $-1,1$ of the interval.

\smallskip

The reconstruction procedure described in \cite{bat} consists of solving
certain linear and non-linear algebraic equations whose coefficients
are expressed through the moments $m_{\alpha}$. These equations have a
very specific structure which we illustrate in the next section with the
simplest example of the classical ``Prony system". The last step requires
also finding a basis of the solution space of the differential equation
$Dg=0$.

\medskip

 We shall apply below this reconstruction procedure,
 referring to it as to Procedure 1.

\subsection{Prony system}

 Prony system appears as we try to solve a very simple version of the shifts
 reconstruction problem.
 Consider $F(x)=\sum_{j=1}^Na_j\delta(x-x_i)$.
 We use as measurements the polynomial moments
 $$
 m_n=\int_{-\infty}^\infty F(x)x^ndx.
 $$
 After substituting $F$ into
  this integral we get
 $m_n=\int \sum_{j=1}^Na_i\delta(x-x_j)x^ndx=\sum_{j=1}^Na_jx_j^n.$
 Considering $a_i$ and $x_i$ as unknowns, we obtain equations
 \be
 m_n=\sum_{j=1}^Na_jx_j^n, \ n=0,1\dots.
 \ee
This infinite set of equations is called Prony system. It can be traced
at least to R. de Prony (1795, \cite{pron}) and it is used in a wide variety
of theoretical and applied fields. See, for example, \cite{bsy,Sar.Yom} and
references therein for a very partial list, as well as for a sketch of one of
the solution methods. This method requires $2N$ equations from (2.2). It allows
first to find the number of nonzero coefficients $a_j$. Then $a_j$ and $x_j$
are found via solving first a Hankel-type
linear system of equations with coefficients formed by the moments $m_l$,
interpreting the solution as the coefficients of a certain polynomial, and
finding all the roots of this polynomial.

 We shall apply below this Prony solution procedure
 in our specific situation, referring to it as to Procedure 2.

\newcommand{\upfun}{\ensuremath{\overline{\phi}}}
\newcommand{\downfun}{\ensuremath{\underline{\phi}}}

\section{Main result}
\setcounter{equation}{0}

We assume that the domain $G \subset  {\mathbb R}^2$ to be reconstructed
has a ``D-finite boundary". More accurately, we have the following definition:

 \bd A compact domain $G \subset  {\mathbb R}^2$ is called $D$-finite if
 the boundary $\partial G$ is a union of $\kappa$ segments
 $S_j$ with the following property: there exists a linear differential
 operator $D$ of the form
 {\rm (2.1)}
 and with the leading coefficient
 not vanishing for $x$ in the projection of $G$, such that each $S_j$
 is the graph of a function $y=\psi_j(x)$ satisfying $D\psi_j=0$.
 \ed
In particular, if each $S_j$ is a graph of an algebraic function
$y=\psi_j(x)$ and all the branches of these algebraic functions are regular
over the projection of $G$, then $G$ is $D$-finite. The simplest but still
important example of $D$-finite domain is when all $\psi_j$ are polynomials.

 We do not restrict the topological type of $G$ --- it may have
 ``holes".

   \medskip

  \begin{figure}
  \centering
  \includegraphics[scale=0.7]{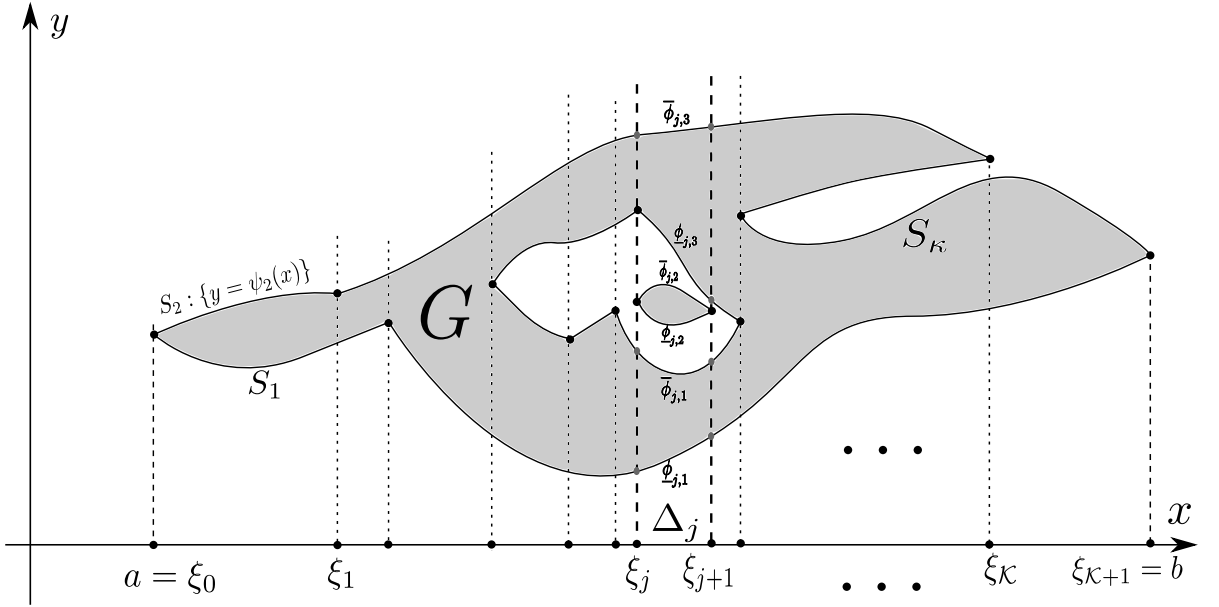}
  \caption{Schematic representation of the domain $G$.}
  \label{fig:fig1}
  \end{figure}

   \medskip

Before we formulate the main result, let us introduce some notations.
Let $[a,b]$ be the projection of $G$ onto the $x$-axis, and let
$a = \xi_0 < \xi_1 < \ldots < \xi_{{\cal K}+1} = b$ be all the projections
of the endpoints of the segments $S_j$ of the boundary $\partial G$
(see Figure 1). Certainly, ${\cal K}\leq \kappa$, while the maximal number
of the intersection points of $\partial G$ with vertical lines is at most
$\kappa-1$.

  Write the linear differential operator $D$ in Definition 3.1 as
 \be
 D g(x) = \sum\limits_{j=0}^{N} p_j(x)\frac{d^j g}{dx^j}(x),
 \ee
  where $p_j(x)= \sum\limits_{i=0}^{k_j} a_{i,j} x^i$.

  \bt Any $D$-finite domain with the discrete parameters
  $\kappa,N,k_j$ can be uniquely reconstructed from a
  collection of double moments
    \be
   m_{\alpha, \beta} = \iint_{R^2} \, \chi_G (x,y)\cdot
   x^{\alpha} y^{\beta} dx dy, \ 0\leq \alpha \leq M(\kappa,N,k_j), \ 0\leq \beta
   \leq 2(\kappa-1).
   \ee
  The reconstruction procedure requires solving of certain
 linear and non-linear algebraic equations whose coefficients
 are expressed through the moments $m_{\alpha,\beta}$, and solving equation
 $ Du=0$
 with specific numerical coefficients found in previous stages.
 \et
 \pr Denote by $\Delta_j$ the interval $[\xi_j,\xi_{j+1}], \ j=0,\dots,{\cal K}.$
 Over each $\Delta_j$ the domain $G$ is a union of $s_j\leq {1\over 2}(\kappa-1)$
 strips $\downfun_{j,l}\leq y \leq \upfun_{j,l}, \ , l=1,\dots,s_j$,
 (see Figure 1). We have
 \be
 m_{\alpha, \beta} = \int_a^b x^\alpha \Psi_\beta(x) dx= \sum_{j=0}^{{\cal K}}
 \int_{\Delta_j}x^\alpha \Psi_{\beta,j}(x)dx,
 \ee
 where for $x\in \Delta_j$
 \be
 \Psi_{\beta,j}(x)=\int_{\downfun_{j,1}(x)}^{\upfun_{j,s_j}(x)} y^\beta \chi_G (x,y)dy=
 {1\over {\beta+1}} \sum_{l=1}^{s_j}
 [\upfun^{\beta+1}_{j,l}(x)-\downfun^{\beta+1}_{j,l}(x)].
  \ee
The first conclusion is that for each $\beta \geq 0$ the function $\Psi_\beta$ is
piecewise $D$-finite. Indeed, on each interval
$\Delta_j, \ j=0,\dots,{\cal K}, \ \Psi_\beta=\Psi_{\beta,j}$ is a linear
combination of the $\beta$-th powers of the functions
$\downfun_{j,l}, \upfun_{j,l}, \ l=1,\dots,s_j$ which, by the assumption on $G$,
satisfy $D\psi = 0$. Hence $\Psi_{\beta}$ (i.e. each $\Psi_{\beta,j}$) satisfies
another linear differential equation with polynomial coefficients
$D_\beta \Psi_\beta=0$. The operator $D_\beta$ depends only on $D$, in particular,
its order and degree depend only on the order and degree of $D$, and it has no
singularities on $[a,b]$, if $D$ possesses this property.

\smallskip

We can find one-dimensional moments of $\Psi_{\beta}$ via (3.3):
 \be
 m_\alpha(\Psi_{\beta})= \int_a^b x^\alpha \Psi_\beta(x) dx=m_{\alpha, \beta}.
 \ee
Now we apply one-dimensional Procedure 1 from Section \ref{sec:picewise-dfin-onedim},
and reconstruct $\Psi_{\beta}$ from
the moments $m_{\alpha, \beta}, \ \alpha =0,1,\dots,\mu_\beta \leq M(\kappa,N,k_j)$
for each $\beta=0,1,\dots,2(\kappa-1)$. The reconstruction procedure requires
solving of certain linear and non-linear algebraic equations whose coefficients
are expressed through the moments $m_{\alpha,\beta}$, and solving differential
equations $D_\beta u=0$ with specific numerical coefficients found in previous
stages. Ultimately, $\Psi_{\beta}(x)=\Psi_{\beta,j}(x)$ is represented on each
interval $\Delta_j$ as a linear combination of the basis solutions of $D_\beta u=0$.

\medskip

Next, for each fixed $x\in \Delta_j$ we consider equalities (3.4) for different
$\beta$ as a system of equations for the unknowns $\downfun_{j,l}(x), \upfun_{j,l}(x)$,
with the known by now right hand side $\Psi_{\beta,j}(x)$. This system of
equations is a special case of the Prony system, as described in the previous
section. Here the amplitudes $a_j$ are known to be $\pm 1$. (3.4) does not give
the first equation in the Prony system, but it is just the sum of the amplitudes $a_j$,
and we know it to be zero. So finally, we apply Procedure 2 and solve system (3.4),
reconstructing the functions $\downfun_{j,l}(x), \upfun_{j,l}(x)$ for each
$x\in \Delta_j$. Now the functions $\downfun_{j,l}(x), \upfun_{j,l}(x)$ for
$j=0,\dots,{\cal K}, \ l=1,\dots,s_j,$ completely determine the domain $G$.
Theorem 3.1 is proved.

\smallskip

 As it was mentioned above, for the case of algebraic boundary segments explicit
 bounds can be given on the number of the moments required for reconstruction. We
 provide these bounds in the case of polynomial boundaries, where the expressions
 are relatively simple and pretty sharp.

 \bt Under the assumptions of Theorem
 {\rm 3.1}
 let us assume additionally that each boundary segment $S_j$
 is a graph of a polynomial $y=\psi_j(x)$ of degree at most
 $d$.
 Then the domain
   $G$
 can be uniquely reconstructed from a collection of double moments
 \[
 \{ m_{\alpha, \beta}:\; 0\leq \beta \leq 2(\kappa-1),\;
 0\leq \alpha \leq M(N,\kappa,\beta) \}
 \]
 where
 \[
 M(d,\kappa,\beta)=
 \max\biggl\{ 2(\beta N + 1)\kappa-2, (\kappa + 1)(\beta N +1) \biggr\}.
 \]
 \et
 \pr The functions $\Psi_\beta$ constructed in the proof of Theorem 3.1
 are now piecewise-polynomials of degree at most
 $\beta N$,
 and the conclusion follows from the corresponding result of
 \cite{bat} given in Section 2 above.

 \section{Examples of Explicit Reconstruction}
 \setcounter{equation}{0}

 The main object considered in this section is a compact plane domain
 bounded by a part of an elliptic curve
 $$
 y^2 = ax^3 + b x^2 + c x + d \equiv f(x).
 \eqno{(4.1)}
 $$

 We assume that the roots $ x_1 < x_2 $, and
 $\, x_3 $ of the equation $ f(x) = 0 $ are real, and that
 the function $f(x) $ is positive on the interval $ (x_1, x_2) $.
 Hence, the equation $ y^2 = ax^3 + b x^2 + c x + d \, $ defines
 a compact domain $ G \subset {\mathbb R}^2$.

 Consider a finite collection of corresponding moments
 $$
 m_{\alpha, \beta} = \int_G x^{\alpha} \cdot y^{\beta} \,dx\, dy.
 $$
 Since the domain $ G $ is symmetric with respect to the
 axis  $ Ox$, we have
 $ m_{\alpha, \beta} = 0 $ for odd values of $ \beta $, and
 $$
 m_{\alpha, 2 \beta} =
 \frac{2}{2 \beta + 1}\int_{x_1}^{x_2} x^{\alpha} \cdot
 (ax^3 + b x^2 + c x + d)^{\beta + 1/2} dx.
 \eqno{(4.2)}
 $$
 In order to simplify the formulae below,
 we introduce one more notation:
 $$
  M_{\alpha, 2 \beta} =
  \frac{2\beta +1}{2} \cdot m_{\alpha, 2 \beta},
  $$
  and we call it "Moment".
 So, the equation (4.2) implies that
 $$
 M_{\alpha, 2\beta +2} = a \cdot M_{\alpha + 3, 2\beta} +
 b \cdot M_{\alpha + 2, 2\beta} + c \cdot M_{\alpha +1 , 2\beta} +
 d \cdot M_{\alpha, 2\beta}.
  \eqno{(4.3)}
 $$

 Our task is to determine the curve (4.1) from a finite
 (possibly minimal) collection of the Moments $ M_{\alpha, 2 \beta} $.

 Let us calculate the Moments
 $ M_{\alpha, 2 \beta} $ for small values of the indices:

 \medskip

  $ M_{0,2} = a \cdot M_{3,0}+ b \cdot M_{2,0} + c \cdot M_{1,0}
  + d \cdot M_{0,0} $; \hfill (4.4)

 $ M_{1,2} = a \cdot M_{4,0}+ b \cdot M_{3,0} + c \cdot M_{2,0} +
 d \cdot M_{1,0} $.

 \medskip

 \noindent
 Hence, one can obtain a system of relations with unknown coefficients
 $ a,b,c,d $, and to verify compatibility of the "data"
  $ \{ M_{\alpha, \beta} \} $.

 \medskip

 Here are two more methods of construction of similar relations:

 a). Consider the Moment
 $$
  M_{\alpha, 2} = \int_{x_1}^{x_2} x^{\alpha} \cdot
  (ax^3 + b x^2 + c x + d)^{3/2} dx .
 $$
  Since $ f(x_1)= f( x_2)=0 $, integrating by parts:

 $ du= x^{\alpha} dx $,
 $ \,\, v= (ax^3 + b x^2 + c x + d)^{3/2}$,
  shows that:
 $$
 M_{\alpha, 2} = -\frac{3}{2(\alpha +1)} \cdot
 (3a\cdot  M_{\alpha+3,0} + 2b\cdot  M_{\alpha+2,0} +
 c \cdot M_{\alpha+1,0}).
 $$
 In the cases $ \alpha = 1$, and $ \alpha = 0$ we get
 $$
   -\frac{4}{3} \cdot M_{1, 2} =
   3a\cdot  M_{4,0} + 2b\cdot  M_{3,0} +
   c \cdot M_{2,0}.
   \eqno{(4.5)}
   $$
 and
 $$
 -\frac{2}{3} \cdot M_{0, 2} =
 3a\cdot  M_{3,0} + 2b\cdot  M_{2,0} +
 c \cdot M_{1,0},
 \eqno{(4.6)}
 $$

 b). Similarly, since $ f(x_1)= f( x_2)=0 $,
 a simple change of the variables implies
 $$
 \int_{x_1}^{x_2} (3 a x^2 + 2 b x + c)
 \cdot (a x^3 + 2 b x^2 + c x + d) dx = 0,
 $$
 and hence,
 $$
 0 = 3 a \cdot M_{2,0} + 2 b \cdot M_{1,0} + c \cdot M_{0,0}.
 \eqno{(4.7)}
 $$

 \medskip

   Consider now the system of three linear equations
  (4.5), (4.6), (4.7) with respect to the unknowns $ a, b, c $.
  Let $\, L_{2,f}[x_1, x_2] \, $ be the Hilbert space composed by
  corresponding functions defined on the segment $ [x_1, x_2] $,
  endowed with scalar product:
  $$
 \langle F(x), H(x) \rangle :=
 \int_{x_1}^{x_2} F(x) \cdot H(x) \sqrt{f(x)} dx.
  $$
 We have in these notations:
 $ M_{4,0} = \langle x^2, x^2 \rangle $,
 $ M_{3,0} = \langle x^2, x^1 \rangle $,
 $ M_{1,0} = \langle x^1, x^0 \rangle $,

 \noindent
 $ M_{2,0} = \langle x^2, x^0 \rangle = \langle x^1, x^1 \rangle $,
 and
  $ M_{0,0} = \langle x^0, x^0 \rangle $.

 \smallskip

 So, the determinant of the system (4.5), (4.6), (4.7) equals
 $$
 6 \cdot \left(
 \begin{array}{ccc}
  M_{4,0} & M_{3,0} & M_{2,0} \\
  M_{3,0} & M_{2,0} & M_{1,0} \\
  M_{2,0} & M_{1,0} & M_{0,0}
  \end{array}
  \right) =
  \, 6 \cdot \left(
 \begin{array}{ccc}
 \langle x^2,x^2\rangle & \langle x^2,x^1\rangle & \langle x^2, x^0\rangle \\
 \langle x^2,x^1\rangle & \langle x^1,x^1\rangle & \langle x^1, x^0\rangle \\
 \langle x^2,x^0\rangle & \langle x^1,x^0\rangle & \langle x^0, x^0\rangle
 \end{array}
  \right),
 $$
 and this coincides (up to the factor 6) with the determinant of the
 Gram matrix of system of three polynomial functions
 $ x^2 $, $ x^1 $, and $ x^0 $, which are linearly independent
 in the space $ L_{2,f}[x_1, x_2] $.
 It is well-known that this determinant is strictly positive, thus,
 the system of linear equations (4.6), (4.7), (4.8) has a unique
 solution $ a,b,c$. Then the coefficient  $\, d \, $ is uniquely
 determined from the equation (4.4), since we assume that $ M_{0,0} > 0 $.

 So, we have proved the following:

 \medskip

 { \bf Theorem 4.1.} (\cite{bgy2})
 {\it In order to reconstruct an elliptic curve
 {\rm (4.1)} it is sufficient to know $7$ moments
 $ m_{0,0} $, $ m_{1,0} $, $ m_{2,0} $, $ m_{3,0} $,
 $ m_{4,0} $, $ m_{0,2} $ è $ m_{1,2} $.
 }
 \medskip

 Note that such an "overdeterminancy" allows to obtain corresponding
 (nonlinear) relations between the moments listed above.

 Similar calculations illustrate theorems 3.1, 3.2 in a very
 simple case:
 ($ {\cal K }= 1 $,
 $ N = 2 $), see \cite{bgy1}.
 Here we reconstructed a triangle
  $ T \subset {\mathbb R}^2 $ with non-vertical edges
  from a given set of moments
 $$ m_{0,0} , \quad m_{1,0} , \quad m_{2,0} , \quad m_{3,0} ,
 \quad  m_{0,1} , \quad m_{1,1} .
  $$
 Cf. also \cite{ghmp}.

 \section{Invisible sets and functions}
 \setcounter{equation}{0}

 Let ${\cal P}_n$ denote the space of polynomials $P(x_1,\dots,x_n)$ and let a
 collection $S\subset {\cal P}_n$ be fixed. We call a function $f$ on
 ${\mathbb R}^n$ $S$-invisible if $\int_{{\mathbb R}^n}P(x)f(x)dx=0$ for each
 $P\in S\subset {\cal P}_n$. A domain $G\subset {\mathbb R}^n$ is $S$-invisible
 together with its characteristic function. With obvious modifications this
 definition is extended to subsets of higher codimension and to distributions.

In this section we discuss some examples of invisible sets and functions, coming
from different fields. Besides Propositions 5.1 and 5.4 below, we do not provide new
results, but rather an initial attempt to find similarity between several
apparently unrelated problems. The common property we would like to stress
is a remarkable ``rigidity" and symmetry of the invisible objects.

\subsection{$S\subset {\cal P}_n$ annihilating a fixed differential operator}

For a fixed partial differential operator ${\cal D}$ in $n$ variables it is
natural to consider $S\subset {\cal P}_n$ consisting of all $P\in S$ with
${\cal D}P=0$.

\subsubsection{Null quadrature domains}

 Put ${\cal D}=\Delta$ to be the Laplacian, and denote $S_h$ the
 corresponding set of harmonic polynomials
 $S_h=\{P \in {\cal P}_n, \ \Delta P=0\}$.

A domain $G\subset {\mathbb R}^n$ is called a null quadrature domain if
$\int_G h dx = 0$ for all harmonic and integrable functions $h$. Taking
$h=P\in S_h$ we get a closely related notion, so null quadrature domains are
essentially all the $S_h$-invisible sets. This class of domains has been
intensively studied, and it has wide applications, in particular, in the
investigation of the Newtonian potential, and of the filtration flow of
incompressible fluid (see \cite{kl1,kl2,Var.Eti} and references therein).

Null quadrature domains include half-spaces, exterior of ellipsoids,
exterior of strips, exterior of elliptic paraboloids and cylinders over
domains of these types. It is known that in ${\mathbb R}^2$ any null
quadrature domain belongs to one of the categories above (\cite{sak}).
A complete description of all null quadrature domains in higher dimensions
has remained an open problem. A significant progress has been recently
achieved (see \cite{kl1,kl2} and references therein).

\subsubsection{Sets invisible for solutions of the wave equation}

Here we consider a somewhat artificial example which however illustrates
the situation for another type of the operator ${\cal D}$. We put $n=2$ and
consider ${\cal D}=W={{\partial^2}\over {\partial x \partial y}}$. In this
case $S_W$ consists of all the polynomials $P$ of the form $P(x,y)=Q(x)+R(y)$.
A function $f(x,y)$ is $S_W$-invisible if and only if
$\int_{{\mathbb R}^2}f(x,y)(Q(x)+R(y))=0$ for any polynomials $Q(x), R(y)$.
In turn, this is equivalent to the vanishing of all the moments

 \smallskip

 $\int_{{\mathbb R}^2} x^k f(x,y)\,dx \, dy=\int x^k dx \int
 f(x,y)dy\, $ and

 $\int_{{\mathbb R}^2} y^l f(x,y)\, dx \, dy =\int y^l dy \int f(x,y)dx.$

  \smallskip

\noindent
  This is equivalent to the identical vanishing of the
  functions

  $F(x)=\int f(x,y)dy=0$ and $H(y)=\int f(x,y)dx=0$. So we
have the following result:

 \bp $f(x,y)\,$ is $\, S_W$-invisible if and only if
 $\, F(x)=\int f(x,y)dy=0\, $ for each $x$ and $\, H(y)=\int f(x,y)dx=0$
 for each $y$. In particular, this is true for functions given by finite
 or infinite sums of the products $\phi(x)\psi(y)$ with
 $\int \phi(x)dx=0, \ \int \psi(y)dy=0.$
 \ep

\subsubsection{Vanishing conjectures of W. Zhao}

In a series of recent papers (\cite{zha1,zha2} and references therein)
W. Zhao has studied a number of vanishing conjectures which relate polynomials
annihilating certain differential operators, invisible sets, and the well
known Jacobian conjecture (\cite{jac}).

For a given ${\cal D}$, specifically, for ${\cal D}=\Delta$
being the Laplacian, the polynomials $P$ have been considered satisfying the
following condition: ${\cal D}^lP^l=0, \ l=1,\dots.$ This condition turned
out to be closely related to the classical and generalized orthogonal polynomials.
The following conjecture has been shown in \cite{zha1} to be equivalent to
the Jacobian conjecture:

\medskip

 {\bf Conjecture A} {\it If for a homogeneous polynomial $P$ of degree four
 $\Delta^lP^l=0, \ l=1,2,\dots$, then $\Delta^lP^{l+1}=0, \ l\gg 1$}.

\medskip

 It was shown in \cite{zha1} that the vanishing of ${\cal D}^lP^l$ is equivalent
 to $P$ being Hessian nilpotent --- i.e. the Hessian matrix
 $H(P)=({{\partial^2 P}\over {\partial x_i \partial x_j}})$ being nilpotent.
 In \cite{zha2} Conjecture A has been closely related to the following
 Conjecture B:

\medskip

 {\bf Conjecture B} {\it For a compact domain $G\subset {\mathbb R}^n$, for
 a positive measure $\mu$ on $G$, and for a polynomial $P$ if all the moments
 $\int_G P^k d\mu$ vanish,
 then
 for any polynomial $q$ the moments
 $\int P^k q \ d\mu$ vanish for $k \gg 1$}.

\medskip

In our language this conjecture can be reformulated as follows: if $(G,\mu)$ is
invisible for all the powers of $P$ then it is ``eventually invisible" for
the sequence of polynomials $P^kq$. Below we discuss this conjecture in somewhat
more detail.

\subsection{Sets invisible by powers of a fixed polynomial}

In this section we discuss the vanishing problem for the moments $\int_G P^k d\mu$,
i.e. the conditions of invisibility of $(G,\mu)$ for all the powers of $P$. Besides
its appearance in Zhao's study of the Jacobian conjecture as above, this question
is related to a wide spectrum of problems in Analysis, Algebra, Differential
Equations, and Signal Processing. We shortly mention below only a very few of these
remarkable connections.

\subsubsection{One-dimensional case}

 In one dimension the question is to describe all the
 univariate polynomials
 (Laurent polynomials, etc)
 $P(x)$ and $q(x)$ for which
 \be
 m_k= \int^b_a P^k(x)q(x)dx = 0, \ k=0,1 \dots .
 \ee
This question appears as a key step in understanding the classical Center-Focus
problem of the Qualitative Theory of ODE's in the case of Abel equation
(see \cite{bfy,broy} and references therein).

Even in this simplest case the answer (only recently obtained in
\cite{mp,pak2}) is far from being straightforward. In particular,
it involves subtle properties of the polynomial composition algebra.
To state the result we need the ``composition condition" (CC) defined
initially in \cite{al} and further investigated in
\cite{bfy,broy,bru.yom,mp,pak1,pak2,pry} and in many other publications.

 \bd Differentiable functions $f(x)$ and $g(x)$ on
 $[a,b] \subset {\mathbb R}$ are said to satisfy a composition
 condition
 {\rm (CC)} on $[a,b]$ if there exists differentiable $W(x)$
 defined on $[a,b]$ with $W(a)=W(b)$, and two differentiable
 functions $\tilde F$ and $\tilde G$ such that $F(x)=\int^x_a
 f(x)dx$ and $G(x)=\int^x_a g(x)dx$ satisfy \be F(x)=\tilde
 F(W(x)), \ G(x)= \tilde G(W(x)), \ x \in [a,b].\ee
 \ed
If $f,g$ are polynomials and they satisfy (CC) then $W$ is necessarily
also a polynomial.

Composition condition implies vanishing of all the moments $m_k$ (change
of variables). Necessary and sufficient condition for vanishing of $m_k$
for $p,q$ polynomials is given by the following theorem:

 \bt {\rm (\cite{mp,pak2})}
 The moments $m_k$ in
 {\rm (5.1)} vanish for
 $k=0,1,\dots$ \textup(i.e. $[a,b]$ is invisible for $P^kq$\textup)
 if and only if \ $q(x)= q_1(x)+\dots+q_l(x),$ with $l=1,2$ or
 $3$, where $q_1,\dots,q_l$ satisfy composition condition
 {\rm (CC)} with $P(x)$
 on $[a,b]$, possibly with different right factors $W_1,\dots,W_l$.
 \et
Analysis of the case of rational functions, and, in particular, of Laurent
polynomials (directly related to the Poincar\'e Center-Focus problem for
plane polynomial vector-fields) turns out to be significantly more difficult
(see \cite{pak1}).

 If we allow $P,q$ above to be only piecewise-polynomial
 (piecewise-rational) then another form of composition condition becomes
 relevant: a ``tree composition condition" (TCC) where $W$ maps $[a,b]$ not
 into ${\mathbb R}$, but into a certain topological tree. Still under some
 restrictions a result similar to Theorem 5.1 remains valid
 (see \cite{bru.yom}). We hope that these recent developments can provide
 a better understanding of invisible $D$-finite domains, as above.

\subsubsection{Some examples in higher dimensions}

 We start with a definition of a multidimensional composition
 condition (MCC) given in \cite{fpyz}, which directly generalizes
 Definition 5.1. (MCC)
 provides a natural {\it sufficient} condition for the moments
 vanishing. However, as we shall see below, in $n>1$ variables this
 condition is much stronger than the vanishing of the ``one-sided"
 moments $m_k= \int_{\Omega} F^k(x) g(x)dx, \ k=0,1,\dots$. In
 fact, it is exactly relevant to the vanishing of the $n$-fold
 moments
 \be
 m_{\alpha}= \int_{\Omega} F_1^{\alpha_1}(x)\cdot ... \cdot F_n^{\alpha_n}(x)g(x)dx,
 \ee
 for all the nonnegative
 multi-indices $\alpha=(\alpha_1,\dots,\alpha_n)$.

\medskip

 Let $\Omega$ be an open relatively compact domain of ${\mathbb R}^n$
 with a smooth boundary $\partial \O$. First we need for maps
 $W: \O \rightarrow {\mathbb R}^n$ a definition generalizing to
 higher dimensions the requirement $W(a)=W(b)$ in dimension one.

 \bd {\rm (\cite{fpyz})}
 A continuous mapping $W: \O \rightarrow {\mathbb R}^n$ is said
 to ``flatten the boundary" $\partial \O$ of $\O$ if the
 topological index of $W|_{\partial \O}$ is zero with respect to
 each point $w\in {\mathbb R}^n \setminus W(\partial \O)$.
 \ed
Informally, $W$ flattens the boundary $\partial \O$ of $\O$ if
$W(\partial \O)$ ``does not have interior" in ${\mathbb R}^n$. In
particular, this is true if $W|_{\partial \O}$ can be factorized
through a contractible $(n-1)$-dimensional space $X$. The simplest
example is when $X$ is a point, so $W$ mapping $\partial \O$ to a
point always flattens the boundary. We have the following simple fact:

 \bp {\rm (\cite{fpyz})}
 A mapping $W: \O \rightarrow {\mathbb R}^n$ flattens the
 boundary $\partial \O$ if and only if the integral
 $\int_\O H(W(x))dW(x)$ vanishes for any function $H(W)$.
 \ep

\medskip

Now let $F_1,\dots,F_s$ be differentiable functions on $\O$ and
let $\mu$ be a measure on $\O$ given by its density $g(x)$:
$d\mu(x)=g(x)dx$.

 \bd {\rm (\cite{fpyz})}
 Functions $F_l, \ l=1,\dots,s$ and a measure $\mu$ on $\O$
 satisfy multi-dimensional composition condition {\rm (MCC)} if there exists
 a differentiable mapping $W: \O \rightarrow {\mathbb R}^n$,
 flattening the boundary $\partial \O$, functions
 $\tilde F_l(w), \ l=1,\dots,s,$ and $\tilde g(w)$ on ${\mathbb R}^n$
 such that $F_l(x)=\tilde F_l(W(x)), \ l=1,\dots,s,$ and
 $\, d\mu(x)=g(x)dx = \tilde g(W(x)) dW.$
 \ed
The following simple proposition (implied directly
by Proposition 5.1) shows that (MCC) is sufficient for moment vanishing:

 \bp If a function $F$ and a measure $\mu$ on $\O$ satisfy
 {\rm (MCC)},
 then
 all the moments $m_k= \int_{\Omega} F^k(x) g(x)dx, \
 k=0,1,\dots,$ vanish.
 \ep
 Consider the following
example: let $\O \subset {\mathbb R}^n$ be defined by $P(x)\leq 1$ for a
certain polynomial $P(x), \ x=(x_1,\dots,x_n) \in {\mathbb R}^n$. For
each $j=1,\dots,n$ define $S_j$ to be a collection of polynomials
$S_j=\{Q^n(P){{\partial P}\over {\partial x_j}}\}$ with $Q$ an arbitrary
univariate polynomial.

 \bp
 For each
 $j=1,\dots,n, \ \ \O, \ Q(P), \ d\mu={{\partial P}\over {\partial x_j}}dx$
 satisfy
 {\rm (MCC)},
 so the domain $\O$ is invisible for $S_j$.
 \ep
 \pr Define $W: \O \rightarrow {\mathbb R}^n$ by
 $W(x_1,\dots,x_n)=(y_1,\dots,y_n)$, with
 $y_i=x_i, \ i\ne j, \ y_j=P(x_1,\dots,x_n)$. $W$ maps the boundary
 $\partial \O$ into the hyperplane $\{y_j=1\}\subset {\mathbb R}^n$, so
 it flattens $\partial \O$. Now, $Q(P)=\tilde Q(W)$, where
 $\tilde Q(y_1,\dots,y_n)=Q(y_j),$ and
  $$
 d\mu={{\partial P}\over {\partial x_j}}dx=dx_1\cdots
 dx_{j-1}\cdot dP \cdot dx_{j+1} \cdots dx_n=dW.
  $$

Let us now describe a situation where (MCC) is a necessary and sufficient
condition for invisibility. Consider double moments of the form
 \be
 m_{k,l}=\int_\O P^k(x,y)Q^l(x,y)r(x,y)dx dy, \ k,l=0,1,\dots, \
 \O\subset {\mathbb R}^2.
 \ee
 We shall assume that $P$ in the domain
 of consideration satisfies ${{\partial P}\over {\partial x}}\ne 0$
 and consider $\O$ of the form $a\leq P(x,y) \leq b, \ c\leq y \leq d$.
 The functions $P,Q,r$ in
 (5.4)
 are assumed to be real analytic, and $Q$ is assumed to have
 a simple critical value on each level curve of $P$ inside $\O$.

 \bt
 {\rm (\cite{fpyz})}
 Under the above assumptions all the moments
 $m_{k,l}, \ k,l=0,1,\dots,$ vanish if and only if $\O,P,Q,r dx$
 satisfy {\rm (MCC)}.
 \et
So the domain $\O$ as above is $S$-invisible for $S$ consisting of
all the products $P^kQ^lr$ if and only if $\O, P,Q,r dx$ satisfy (MCC).

\subsubsection{Mathieu conjecture and Laurent polynomials}

Zhao's Conjecture B above has been motivated, in particular, by the
following conjecture of O. Mathieu (\cite{mat}), closely related to
many important questions in Representation Theory: let $M$ be a
compact Lie group. Denote $F(M)$ the set of $M$-finite functions
on $M$ (i.e. polynomials in all the characters on $M$) and let
$\mu$ be the Haar measure on $M$.

\medskip

\noindent{\bf Conjecture C}. {\it If for some $f(x)\in F(M)$
 \be
 \int_M f^k(x)d\mu(x) = 0, \ k=1,2,\dots
  \ee
  then for any $g(x)\in F(M)$ we have
  $\int_M f^k(x)g(x)d\mu(x) = 0, \ k \gg 1$}.

\smallskip

This conjecture is known to imply the Jacobian conjecture (\cite{mat}).
In our language, it states that if $M$ is invisible for $f^k$, it is
eventually invisible for $f^kg$ with any $M$-finite function $g$.

Conjecture C has been verified in \cite{dui} for the Abelian
$M$, i.e. for $M$ being the $n$-dimensional torus $T^n$. In this
case $M$-finite functions are Laurent polynomials in
$z=(z_1,\dots,z_n), \ z_i \in {\mathbb C}, \ \vert z_i \vert =1$.
In fact, the following result has been established in \cite{dui}:

 \bt Let $f(z_1,\dots,z_n)$ be a Laurent polynomial. Then the constant
 term of $f^k$ vanishes for $k=1,2,\dots$ if and only if the convex
 hull of the support of $f$ does not contain zero.
 \et
 Here the support of $f$ is the set of multi-indices of all the monomials
 in $f$ with nonzero coefficients. Theorem 5.3 immediately implies
 Conjecture C since under its conditions the support of $f^k$ eventually
 gets out of any compact set on ${\mathbb Z}^n$, in particular, out of
 the support of $g$.

\smallskip

Recently a rather accurate description of moment vanishing
conditions for one-dimensional rational functions and, specifically,
for Laurent polynomials has been obtained in \cite{pak1}. In particular,
an extension of the result of Duistermaat and van der Kallen
(\cite{dui}, Theorem 2.1 above) obtained in \cite{pak1} provides
such conditions:

 \bt
 {\rm (\cite{pak1}, Theorem 6.1)}.
 Let $L(z)$ and $m(z)$ be Laurent
 polynomials such that the coefficient of the term $1\over z$ in
 $m(z)$ is distinct from zero. Assume that $\int_{S^1} L^k(z)m(z)dz
 =0, \ k\gg 1$. Then $L(z)$ is either a polynomial with zero constant
 term in $z$, or a polynomial with zero constant term in $1\over z$.
 \et
As it was explained above, this property implies that
$\int_{S^1} L^k(z)h(z)dz =0, \ k\gg 1$ for any Laurent polynomial $h(z)$.
In particular, we get $$\int_{S^1} L^k(z)g(z)m(z)dz =0, \ k\gg 1$$ for
any Laurent polynomial $g(z)$. Therefore Zhao's Conjecture A holds for
$S^1$ and the measure $d\mu(z)=m(z)dz$.

\smallskip

In \cite{pry} under a stronger assumption of vanishing of the moments
starting from the initial indices, we get the same conclusion assuming
that only a ``horizontal strip" of the moments vanish.

\subsection{Complex moments}

Problems of reconstruction of sets and functions from their complex
moments, and, in particular, vanishing conditions for complex
moments form an important field of investigation in Several Complex
Variables, in Inverse Problems in PDE's and in related fields. We give
here only a few examples illustrating connections with our setting.

\subsubsection{Wermer's theorem and later developments}

The classical theorem of Wermer (\cite{wer}) gives conditions for vanishing
of all complex moments $\int_\gamma x^iy^jdx$ for a closed curve
$\gamma \subset {\mathbb C}^2$: this happens if and only if $\gamma$ bounds
a compact complex one-chain. See \cite{dh,hl1,hl2,wer} and references therein
for an accurate statement and further developments. In our terms $\gamma$ is
invisible for all complex moments if and only if it bounds a compact complex
one-chain.

The theorem of Dolbeault-Henkin (\cite{dh}) gives a remarkable extension of
Wermer's condition to the case of curves $\gamma$ bounding a compact complex
one-chain in the projective space. In particular, in such case the moment
generating function satisfies a non-linear Burgers-type partial differential
equation. This last fact can be reinterpreted as an invisibility of $\gamma$
for certain combinations of the complex moments.

\medskip

Let's assume now that $\gamma \subset {\mathbb C}^2$ is an image of a not
necessarily closed curve $\sigma  \subset {\mathbb C}$ under a rational mapping
$(P,Q):\sigma \rightarrow {\mathbb C}^2$. In this case a more accurate form of
Wermer's theorem can be obtained:

 \bt
 {\rm (\cite{pak1}, Theorem 5.2)}
 The moments $m_{i,j}=\int_\gamma x^iy^jdx$ vanish
 for $i,j\gg 1$ if and only if there exist rational functions $\tilde P, \tilde Q, W$
 such that $P(z)=\tilde P(W(z)), Q(z)=\tilde Q(W(z))$, the curve
 $\bar \sigma=W(\sigma)\subset {\mathbb C}$ is closed, and all the poles of
 $\tilde P, \tilde Q$ lie on one side of the curve $\bar \sigma$.

 In particular, if the moments $m_{i,j}$ vanish for $i,j\gg 1$, they
 in fact vanish for all $i,j>0$.
 \et
In our terms, if the curve $\gamma =(P,Q)(\sigma)$ is eventually invisible for
the complex moments $m_{i,j}$ then it is closed, it is an image under
$\tilde P, \tilde Q$ of the closed curve $\sigma \subset {\mathbb C}$, and
$\gamma$ bounds the
compact complex one-chain $(\tilde P, \tilde Q)(G)$ where $G\subset {\mathbb C}$
is the domain bounded by $\sigma$ and free of poles of $\tilde P, \tilde Q$.

\subsubsection{Complex moments of planar domains}

This problem has been intensively studied, in particular, in relation with
the filtration flow of incompressible fluid, and with the inverse problem
of two-dimensional Potential Theory (see \cite{Var.Eti,kl1,kl2} and references
therein). In particular, in \cite{Var.Eti} one can find a discussion of the
non-uniqueness of reconstruction.

For reconstruction of polygonal domains see \cite{ggmpv,ghmp}. In general,
an important class of quadrature domains (see \cite{as,gp} and references
therein) provides a natural framework for a study of the ``finite dimensional"
reconstruction problem, as well as its possible non-uniqueness, in particular,
the invisibility phenomenon.

 \section{Conclusions}
 \setcounter{equation}{0}

The results of the present paper leave open some important
questions:

\smallskip

 1. We have insisted on a requirement that the singular points of the
 differential operator $D$ are outside the projection of the domain
 to be reconstructed. This requirement excludes some important
 classes of domains. In particular, it prevents reconstruction of
 general semi-algebraic domains $G$ on the plane. Indeed, typically
 the projection of $G$ to the $x$-axis will have singularities on the
 boundary $\partial G$. These singularities will be necessarily also
 singularities of the differential operator $D$ annihilating the
 corresponding algebraic functions. The construction of \cite{bat}
 allows for an adaptation to singular situations: just, on the last
 step we have to use the bases of ``singular bounded solutions" of
 $D=0$. See, for example, \cite{ABK}.

\smallskip

2. An important parameter entering the procedure of reconstruction
of a $D$-finite function $f$ is the maximal number $\mu$ of its
initial moments that can vanish, unless all the moments vanish
identically (see an example in Section 2.1 above). Recently we've
shown that this number depends only on the combinatorial data in
case where singular points of $D$ differ from the jump-points of
$f$. If some of the singularities of $D$ are the jump-points, we
expect an explicit bound on $\mu$ through the size of the
coefficients of $D$.

\smallskip

3. Robustness estimates. A ``quantitative version" of the question
of bounding of the number $\mu$ is to bound all the moments of $f$
through its first $\mu$ moments (``moments domination''). We expect
this problem to be a central one for the robustness estimates of the
reconstruction procedure. Indeed, moments domination implies
directly a bound on the norm of $f$ itself through its first $\mu$
moments. This question is also directly related to the analysis of
periodic solutions of Abel equation (\cite{broy}).

\vskip1cm


\bibliographystyle{amsplain}

\end{document}